\documentclass[preprint,12pt,article,english]{article}

\usepackage{amssymb}
\usepackage{amsmath}

\def\ol{\overline}

\def\s{\mathbf{s}}

\newcommand{\cl}{\mathrm{cl}\,}

\def\R{\mathbb{R}}
\def\N{\mathbb{N}}
\def\eps{\varepsilon}

\def\supp{\mathrm{supp}\,}
\def\qed{$\square$}

\newtheorem{thm}{Theorem}[section]
 \newtheorem{definition}[thm]{Definition}
 \newtheorem{lemma}[thm]{Lemma}
  \newtheorem{corollary}[thm]{Corollary}
   \newtheorem{proposition}[thm]{Proposition}

\begin{document}

\title{Inverse Function Theorem in Fr\'echet Spaces\footnote{
Research is supported by
 the Bulgarian National Scientific Fund under grant KP-06-H22/4.
}}


\author{Milen Ivanov\footnote{Radiant Life Technologies Ltd., Nicosia, Cyprus, e-mail: milen@radiant-life-technologies.com} \and Nadia Zlateva\footnote{Faculty of Mathematics and Informatics,               St. Kliment Ohridski University of Sofia,
              5, James Bourchier blvd., 1164 Sofia, Bulgaria, e-mail: zlateva@fmi.uni-sofia.bg}}

\maketitle

\begin{abstract}
We consider the classical Inverse Function Theorem of Nash and Moser from the angle of some recent development by Ekeland and the authors.

Geometrisation of tame estimates coupled with certain ideas coming from Variational Analysis when applied to a directionally differentiable function, produce very general surjectivity result and, if injectivity can be ensured, Inverse Function Theorem with the expected Lipschitz-like continuity of the inverse.

We also present a brief application to differential equations.
\end{abstract}



\section{Introduction}\label{sec:introduction}

Theorem of Nash and Moser is a powerful tool for studying nonlinear problems with infinitely smooth data. Since the spaces of infinitely smooth functions are not Banach, the usual Inverse Theorems do not work. Good overview of the subject is the monograph \cite{pdonmt}.

The Inverse Function Theorem of Nash and Moser is tough to be proved as it can be seen from the notorious survey \cite{hamilton}. Its proof relies on Newton method and in order for it to work the function should be  smooth and the spaces \textit{tame}.

Ekeland \cite{ekeland} has proved surjectivity result for functions which are only G\^ateaux differentiable. The results of \cite{ekeland} were extended to the case of non-autonomous tameness estimates in \cite{esere}. Since the method of \cite{ekeland} does not require second derivatives, it can be extended to multi-valued maps, see \cite{ngth}. We further this development by proving surjectivity result for multi-valued maps with estimates of all seminorms in \cite{NME-JOTA}. A simple and different proof in the case of merely directionally
differentiable function in  Fr\'echet-Montel spaces is given in \cite{iz-dban}. Recently Inverse Function Theorem in Fr\'echet-Montel spaces was proved in \cite{fabian-cibulka}.

Here we combine ideas from \cite{NME-JOTA} and \cite{iz-dban} to prove an Inverse Function Theorem.

This theorem is weaker than Nash-Moser Theorem, see \cite{hamilton}, because the \textit{tameness} estimates do not depend on the variable. On the other hand, however, the regularity condition we impose on the image Fr\'echet space is very mild.

\begin{thm}\label{thm:main-main}
 Let $(X,\|\cdot\|_n)$ be a Fr\'echet space and $(Y,\|\cdot\|_n)$ be a non-exotic Fr\'echet space and let
  $
    f:X\to Y
  $
  be a continuous function. Let  $U$ be a nonempty  open subset of $X$ such that $f$ is injective and  directionally differentiable on $U$; and there are $d\in \N$ and $c_n>0$ such that for any $x\in U$ and any $v\in Y$,
  \begin{equation}\label{oh}
  \exists h\in X\ :\ f'(x,h)=v\ \mbox{and}\ \| h\|_n\le c_n\|v\|_{n+d},\ \forall n\ge 0.
  \end{equation}
  Then for each $x\in U$ there exists an open $V\ni x$ such that $f$ is invertible on $V$ and
  \begin{equation}
     \label{eq:invert}
     \|f^{-1}(u)-f^{-1}(v)\|_n\le c_n\|u-v\|_{n+d},\quad \forall u,v\in f(V),\ \forall n\ge 0.
  \end{equation}
   \end{thm}

Let us note that the presented here results do not follow from the results for multi-valued maps in \cite{NME-JOTA} and are more general than the results in \cite{iz-dban}.

The paper is organized as follows. In Section~\ref{sec:preliminaries} we recall the necessary definitions and prove some auxiliary results. In Section~\ref{sec:surjectivity} we prove very general surjectivity result. However, in this setting the estimates corresponding to \eqref{eq:invert} hold only if closure is taken. In Section~\ref{sec:main} we prove our main result -- Theorem~\ref{thm:main-main}. In Section~\ref{sec:injectivity} we give some conditions in the spirit of \cite{hamilton} ensuring injectivity, but more research is this direction is required even in Banach spaces case.

   We do not consider here the natural further question of global invertibility.

   In Section~\ref{sec:ift} we follow the standard route of deriving Implicit Map Theorem from uniform surjectivity and then, in the final Section~\ref{sec:ode}, we apply it to generalise an application found in  \cite{ngth}, to a Cauchy problem.

\section{Preliminaries}\label{sec:preliminaries}

Generally speaking, a \textit{Fr\'echet space} is a complete locally convex topological vector space whose topology can be generated by a translation-invariant metric. Then taking Minkowski functions of a countable local base of convex symmetric neighbourhoods of zero, we obtain countable family of seminorms, which also define the topology, see for details \cite[pp.110-114]{fabetal}.

In the context of Nash-Moser-Ekeland theory, however, a set of seminorms is given in advance and the estimates are \textit{in terms} of the given seminorms. Therefore, when the seminorms are fixed, we simply say that \textit{Fr\'echet space} $(X,\|\cdot\|_n)$ is a linear space $X$ with a collection
of seminorms $\|\cdot\|_{n}$, $n=0,\ldots,\infty$, which is separating, that is,
$\|x\|_{n}=0$, $\forall n,$ if and only if $x=0$; and, moreover, equipped with  the metric
\begin{equation}\label{eq:rho}
\rho_{(X,\|\cdot\|_n)}(x,y):=\max_{n\ge 0}\frac{2^{-n}\|x-y\|_{n}}{1+\|x-y\|_{n}}
\end{equation}
$(X,\rho_{(X,\|\cdot\|_n)})$ is complete metric space.

The closed ball centered at $x$ with radius $r$ in $(X,\rho_{(X,\|\cdot\|_n)})$ is denoted by
\[ B(x,r):=\left\{ y\in X : \rho_{(X,\|\cdot\|_n)}(x,y)\le r\right\},\]
and $B^\circ(x,r)$ is the open ball.

\begin{definition}\label{dfn:non-exotic}
 The Fr\'echet space $(X,\|\cdot\|_n)$ is called non-exotic if dropping finitely many seminorm does not change the topology. In other words, for any $N\in\mathbb{N}$ the metric $\rho_{(X,\|\cdot\|_n)}$, see \eqref{eq:rho}, and the metric
 $$
     \rho_{N(X,\|\cdot\|_n)}(x,y)=\max_{n\ge N}\frac{2^{-n}\|x-y\|_{n}}{1+\|x-y\|_{n}}
 $$
 produce the same convergence.
\end{definition}

It is clear that the examples of  Fr\'echet space which do not satisfy Definition~\ref{dfn:non-exotic} would be very \textit{exotic}. The above definition does not in any way exclude Banach spaces, since in that case all seminorms can be taken equal to the original norm, so dropping few of them changes nothing. Also, in many applications some monotonicity of the seminorms (like, e.g. $\|\cdot\|_n\le \|\cdot\|_{n+1}$ for all $n$) is assumed, and such so called \emph{graded spaces} are clearly non-exotic. In short, we need this condition only for accommodating the so called \textit{loss of derivatives}, and it is not that restrictive.

Perhaps, the most used example of non-exotic (and non-Banach) Fr\'echet space is $X = C^\infty (\Omega)$, where $\Omega$ is a compact domain in $\mathbb{R}^n$.

Denote (as somewhat standard)
$$
  \mathbb{R}_+^\infty := \{(s_{n})_{n=0}^{\infty}:\ s_{n}\ge0\},
$$
that is, $\mathbb{R}_+^\infty$ is the cone of all positive sequences indexed from $0$ on.

  For $\mathbf{s}\in \mathbb{R}_+^\infty$ set
\[    \supp \mathbf{s} : = \{n\ge 0:\ s_n > 0\},\]
 \[   |\mathbf{s}| : = \max_{n\ge 0} \frac{2^{-n}s_{n}}{1+s_{n}}.\]
  For a Fr\'echet space $(X,\|\cdot\|_n)$ and a given $\mathbf{s}\in \mathbb{R}_+^\infty$ also define
\begin{equation}
	\label{eq:def:pi}
	\Pi_{\mathbf{s}}(X):=\{x\in X:\ \|x\|_{n}\le s_{n},\ \forall n\ge 0\},
\end{equation}
and
  \begin{equation}\label{eq:def:X:s}
    X_{\mathbf{s}}:=\bigcup_{t\ge 0}t\Pi_{\mathbf{s}}(X),
  \end{equation}
see \cite[Definition 3.1]{NME-JOTA}.

We will now discuss some properties of the while ago introduced  structure  $\Pi_{\mathbf{s}}(X)$ that will be used in the sequel.

It is easy to check that $\Pi_\s(X)$ is closed in $X$ because of
\[
\lim _{k\to \infty}\rho_{(X,\|\cdot\|_n)}(x_k,x)=0 \Leftrightarrow \lim_{k\to \infty} \| x_k-x\|_n =0,\ \forall n\in \N.
\]

What is obvious is that, if $x\in X$ and $\mathbf{s}=(\|x\|_n)_{n=0}^{\infty}$, then
\[  |\mathbf{s}| = \rho_{(X,\|\cdot\|_n)}(0,x),\]
where $\rho_{(X,\|\cdot\|_n)}$ is defined by \eqref{eq:rho}.

The following is one of the key relations in the approach we present.
\begin{lemma}\label{lem:key}
  If $(X,\|\cdot\|_n)$ is Fr\'echet space and $\rho_{(X,\|\cdot\|_n)}$ is defined as in \eqref{eq:rho}, then
  $$
    c\Pi_{\mathbf{s}}(X)\subset B(0, c |\mathbf{s}|)
  $$
  for any $\mathbf{s}\in \mathbb{R}_+^\infty$ and any $c\ge 1$.
\end{lemma}
 For a proof, see \cite[Lemma 3.2]{NME-JOTA}.

Let $(X,\|\cdot\|_n)$ be a Fr\'echet space, let $\mathbf{s}\in \mathbb{R}_+^\infty$ be such that
$\supp \mathbf{s}\ne\varnothing$ and let $X_\mathbf{s}$ be defined by $(\ref{eq:def:X:s})$.
  For $x\in X_\mathbf{s}$   define
 \[   \|x\|_{\mathbf{s}}:=\sup \left\{\frac{\|x\|_{n}}{s_{n}}:\ n\in\supp \s\right\}.\]
  Then, $\Pi_{\mathbf{s}}(X)$ is the unit ball of the norm $\|\cdot \|_{\mathbf{s}}$ and $(X_{\mathbf{s}},\|\cdot\|_{\mathbf{s}})$ is a Banach space,
see \cite[Lemma 3.4]{NME-JOTA}.

\begin{proposition}\label{pro:oct_1}
  Let  $(X,\|\cdot\|_n)$ and $(Y,\|\cdot\|_n)$ be Fr\'echet spaces and let $U\subset X$  be nonempty and open.

  Let the function $f:X\to Y$ be continuous and  such that
\begin{equation}
\label{eq:oct_2}
 \cl f(x+\Pi_{\mathbf{s}}(X))\supset f(x)+ \Pi_{ \mathbf{s}}(Y),\quad\forall x:\ x+\Pi_{\mathbf{s}}(X)\subset U
\end{equation}
  Then
    \[
  f (B(x,r))\supset B^\circ (f(x),r),\quad
  \forall x\in U,\ \forall r:\ 0<r< m_U(x),
  \]
  where
  $$
  m_U(x) := \mathrm{dist}(x,X\setminus U).
  $$
\end{proposition}

\textbf{Proof. }
 We refer to the proof of Theorem~2.2 from \cite{NME-JOTA}.

 There $V=Y$ and the graph of $f$ is indeed closed, because $f$ is continuous. The condition \eqref{eq:oct_2} means that, in terms of \cite{NME-JOTA}, $f$ is weakly $\Pi$-surjective with $\kappa = 1$.

 Fix $x\in U$ and $r>0$ such that $r < m_U(x)$.

 The proof of \cite[Theorem~2.2]{NME-JOTA} says that for each $\alpha\in (0,1)$
 $$
    \cl f(B(x,r))\supset B(f(x),\alpha r).
 $$
 Then \cite[Lemma~3.5]{NME-JOTA}, or \cite[Theorem~2.55]{Ioffe-book}, gives that
$$
     f(B(x,r))\supset B(f(x),\beta r), \ \forall \beta \in (0,\alpha).
 $$
 But, since $\alpha \in (0,1)$ was arbitrary, the latter in effect means that $f(B(x,r))$ contains $B(f(x),\beta r)$ for all $\beta \in (0,1)$, or in other words,
$$
f (B(x,r))\supset B^\circ (f(x),r).
$$
\qed

\begin{definition}
  Let $(X,\rho_X)$, $(Y,\rho_Y)$, be linear metric spaces, $U\subset X$ be an open set, and
  $$
   f:X\to Y
  $$
  be a function. The  derivative of $f$ at the point $x\in U$ in the direction $h\in X$ is defined by
  $$
    f'(x,h):=\rho_Y{-}\lim _{t\downarrow 0} \frac{f(x+th)-f(x)}{t} .
  $$
  We say that $f$ is directionally differentiable at $x\in U$ if the directional derivative of $f$ at $x$ exists in any direction, and that  $f$ is directionally differentiable on $U$ if it is directionally differentiable at any $x\in U$.

  If $f$ is directionally differentiable on $U$  and the  derivative $f'(x,h)$ is continuous jointly in $(x,h)\in U\times X$, then $f$ is said to be smooth on $U$, denoted by $f\in C^1(U)$.
\end{definition}

If $f\in C^1(U)$ then necessarily $f'(x,h)$ is linear in $h$, denoted by $f'(x,h)=f'(x)(h)$, see \cite[Theorem 3.2.5]{hamilton}.

In \cite{LOEV} we have established the so called \emph{Long Orbit or Empty  Value (LOEV) Principle}   and we have used it there for getting surjectivity results in Banach spaces, and in \cite{NME-JOTA} in Fr\'echet spaces. Here we will recall it and use it as well.

Let $\left(M,\rho\right)$ be a complete metric space.

  Let $S:M\rightrightarrows M$ be a multi-valued map. We say that $S$
  \textbf{satisfies the condition} $\left(\ast\right)$ if $x\notin S(x)$, $\forall x\in M$,
  and whenever $y\in S(x)$ and $\lim_{n}x_{n}=x$, there are infinitely
  many $x_{n}$'s such that
  $
  y\in S(x_{n}).
  $

Informally, LOEV Principle  says that under $\left(\ast\right)$ there exist long $S$-orbits or points where $S$ is empty-valued, thus the name. We will use the LOEV Principle here in the following form
\begin{proposition}[Corollary 3, \cite{LOEV}]\label{thm:loev}
Let $S:M\rightrightarrows M$ satisfy $\left(\ast\right)$
and let $x_{0}\in M$ and $K > 0$ be arbitrary.

Then at least one of the two conditions below is true:

(a) There are $x_{i}\in M$, $i=1,2,\ldots,n+1$, such that
\[
x_{i+1}\in S(x_{i}),\quad i=0,1,\ldots,n,
\]
and
\[
\sum_{i=0}^{n}\rho(x_{i},x_{i+1})>K;
\]

(b) There is $\ol{x}\in M$ such that $\rho(x_{0},\bar{x})\le K$
and
\[
S(\ol{x})=\varnothing.
\]
  \end{proposition}

  We finish this section of miscellaneous auxiliary results by a well-known mean value inequality for  upper Dini derivative.

  \begin{lemma}
   \label{lem:dini-mvt}
   Let $g:[0,1]\to \mathbb{R}$ be a continuous function such that $g(0)=0$ and
   $$
     g^+(\lambda):=\limsup_{t\downarrow 0} \frac{f(\lambda+t)-f(\lambda)}{t}\le 1,\quad \forall \lambda\in(0,1).
   $$
   Then, $g(1)\le 1$.
  \end{lemma}

\textbf{Proof.} 
  Let $a>1$ be arbitrary and let
  $$
    I := \{ x\in [0,1]:\ g(s)\le as,\forall s\in [0,x]\}.
  $$
  Because $g$ is continuous, $I$ is a closed interval, say $[0,b]$. If $b<1$ then for any $t>0$ small enough $b+t<1$ and for any such $t>0$ there exists $t'\in (0,t)$ such that $g(b+t') > a(b+t')$. So $g(b+t')-g(b)> a(b+t') - ab = at'$. Therefore, $g^+(b)\ge a$, contradiction. Hence, $b=1$ and $g(1)\le a$.
  \qed

\section{Local surjectivity}\label{sec:surjectivity}

Here we prove a general surjectivity result. It is established through combining ideas from \cite{NME-JOTA} and \cite{iz-dban}.

If $(X,\|\cdot\|_n)$ is a Fr\'echet space, define
\begin{equation}\label{eq:norm-k}
    |\cdot| _k:=\max_{0\le i\le k}\| \cdot \|_i.
\end{equation}
 Then $(X, |\cdot |_k)$ is a graded Fr\'echet space. We will also use
 \begin{equation}\label{eq:rho-k}
   \rho_{(X,\|\cdot\|_n)}^k (x,y) := \max_{0\le n\le k}\frac{2^{-n}\|x-y\|_{n}}{1+\|x-y\|_{n}}.
 \end{equation}
It is clear that $\rho_{(X,\|\cdot\|_n)}^k $ tends uniformly to $\rho_{(X,\|\cdot\|_n)}$, see \eqref{eq:rho}, in the sense that for any $\varepsilon >0$ there is $N\in \N$ such that
$$
   \left|\rho_{(X,\|\cdot\|_n)}^k (x,y)-\rho_{(X,\|\cdot\|_n)}(x,y)\right|<\varepsilon, \quad\forall k>N,\ \forall x,y\in X.
$$

\begin{lemma}\label{lem1}
Let $(X,\|\cdot\|_n)$ be a Fr\'echet space, and let $A\subset X$ be a non-empty set. Let $x\in X$ be such that
\begin{equation}\label{eq:lem1}
    \forall \varepsilon>0,\forall k\ge 0,\ \exists a\in A:\quad |x-a|_k<\varepsilon.
\end{equation}
Then $x\in \cl A$.
\end{lemma}

\textbf{Proof.} 
Fix $\varepsilon>0$. Fix $k\ge 0$ so large that $\left|\rho_{(X,\|\cdot\|_n)}^k -\rho_{(X,\|\cdot\|_n)}\right|<\varepsilon$ uniformly. Take $a\in A$ such that $|x-a|_k<\varepsilon$. It is clear from \eqref{eq:norm-k} and \eqref{eq:rho-k} that $\rho_{(X,\|\cdot\|_n)}^k (x,a)\le |x-a|_k$. So, $\rho_{(X,\|\cdot\|_n)}(x,a) < \rho_{(X,\|\cdot\|_n)}^k (x,a) + \varepsilon < 2\varepsilon$.

Since $\varepsilon>0$ was arbitrary, $x$ is at zero distance from $A$, or  which is the same,  $x$ is in its closure.
\qed

It is easy to check that the ``only if'' direction is also true, but we are not going to use it.

In a Banach space $(X,\|\cdot\|)$ we denote the closed unit ball by
\[
B_X:=\{x\in X:\|x\|\le 1\}.
\]

\begin{proposition}\label{prop:4}
Let $(X,\|\cdot\|)$ be a Banach space and let $(Y,\|\cdot\|_n)$
be a Fr\'echet space.
Let $g:X\to Y$ be a continuous function which is directionally differentiable  and such that $g(0)=0$.

  If for some $ \ol y\in Y$
 there exist $\sigma >0$ and an open set $U\supset \sigma B_X$   such that
\[
\ol y\in \{g'(x,h):\ h\in\sigma B_X\},\quad\forall x\in U,
\]
then
\[
\ol y\in\cl {g(\sigma B_X)}.
\]
\end{proposition}

\textbf{Proof.} 
Obviously, we may assume that $\ol y\ne 0$.

We will apply Lemma~\ref{lem1} to the set $g(\sigma B_X)$ in $Y$. We consider $Y$ with the graded seminorms $|\cdot|_k$ as defined in \eqref{eq:norm-k} and equip it by the metric $\rho_{(Y,|\cdot|_k)}$ defined as in \eqref{eq:rho} by using these seminorms.

Fix arbitrary $\varepsilon_1>0$.

Let $k_0$ be such that $| \ol y |_{k_0}\neq 0$. Fix arbitrary $k\ge k_0$, $k\in \N$.

Take $\varepsilon\in(0,|\ol y|_k)$ such that
$$
    \varepsilon(1+|\ol y|_k)<\varepsilon_1.
$$
Fix $\mu$ such that
\begin{equation}
  \label{eq:mu-def}
  \mu> \sigma > (1 - \varepsilon)\mu.
\end{equation}

 We intend to apply Proposition~\ref{thm:loev} to $\cl U$, which is complete in the metric induced by $\|\cdot\|$, with $x_0=0$ and $K=\sigma$ for the map  $S:\cl U\rightrightarrows U$ defined for $x\in U$ by
 \begin{equation}
  \label{eq:S-def}
    S(x) := \{u\in U:\ \exists t\in (\mu^{-1}\|u-x\|,\varepsilon):\  |g(u)-g(x)-t\ol y|_k<\varepsilon t\},
 \end{equation}
and $S(x):=\varnothing$ for $x\not\in U$.

We have to show that $S$ satisfies $(\ast)$. Indeed, if $u=x$ in \eqref{eq:S-def} then $t|\ol y|_k < \varepsilon t$, contradiction with the choice of $\varepsilon$. Thus, $x\not\in S(x)$, $\forall x\in\cl U$.

The other requirement of $(\ast)$ follows by continuity:
if $u\in S(x)$ and $x_n\to x$ then for  $t$ corresponding
to $u$ in the definition of $S(x)$, see \eqref{eq:S-def}, we will have $t\in (\mu^{-1}\|u-x_n\|,\varepsilon)$ and
$|g(u)-g(x_n)-t\ol y|_k<\varepsilon t$
for $n$ large enough.
Thus $S$ satisfies $(\ast)$.

Moreover,
\begin{equation}
  \label{eq:S-no-empty}
  S(x)\neq \varnothing,\quad\forall x\in U.
\end{equation}

Indeed, fix $x\in U$. By assumption, there is $h\in X$ with $\|h\|\le\sigma$ such that $g'(x,h) = \ol y$. That is,
\[
\lim_{s\downarrow 0}\rho_Y\left(\frac{g(x+sh)-g(x)}{s},\ol y\right)= 0,
\]
hence
\[
\lim_{s\downarrow 0}\left|\frac{g(x+sh)-g(x)}{s}-\ol y\right|_k= 0.
\]
So, for small $s$, we have that $s\in (0,\eps) $, $x+sh\in U$ and
\[
\left|g(x+sh)-g(x)-s\ol y\right|_k<s\eps .
\]

For $u:=x+sh$ we have that
\[
\frac{\| u-x\|}{\mu}=\frac{s\| h\|}{\mu}\le \frac{s\sigma}{\mu}<s,
\]
and
\[
|g(u)-g(x)-s\ol y|_k<\eps s.
\]
Hence, $u\in S(x)$ and \eqref{eq:S-no-empty} is verified.

Since $\sigma B_X\subset U$, \eqref{eq:S-no-empty} means that (b) from Proposition~\ref{thm:loev} is not an option, and, therefore, there is an $S$-orbit $x_0=0$, $x_1,\ldots,x_{n+1}$, $x_{i+1}\in S(x_i)$, for  $i=0,1,\dots,n$, such that
\begin{equation}
   \label{eq:long-orbit}
   \sum_{i=0}^n \|x_{i+1}-x_i\| > \sigma.
\end{equation}
Denote by $t_i$ some  $t$ from the definition of $x_{i+1}\in S(x_i)$, see \eqref{eq:S-def}, so
\begin{equation}
   \label{eq:t-ineq}
   \|x_{i+1}-x_i\| < \mu t_i,\quad |g(x_{i+1})-g(x_i)-t_i\ol y|_k <\varepsilon t_i.
\end{equation}
For each $i=0,1,\ldots n$  set for brevity
$$
  p_i := \sum_{j=0}^it_j .
$$
From \eqref{eq:t-ineq} we get
\begin{eqnarray*}
\left|g(x_{i+1})-p_i\ol y\right|_k&=&\left|\sum_{j=0}^{i}(g(x_{j+1})-g(x_{j}))- \sum_{j=0}^{i}t_j\ol y\right|_k\\
&\le & \sum_{j=0}^{i} \left|g(x_{j+1})-g(x_{j})- t_{j}\ol y\right|_k \\
 & < &  \sum_{j=0}^{i}\varepsilon t_j=\varepsilon p_i.
\end{eqnarray*}
So,
\begin{equation}\label{eq:y:smaller}
 \left|g(x_{i+1})-p_i\ol y\right| _k< \varepsilon p_i,\quad\forall i=0,1,\ldots,n.
\end{equation}

On the other hand, from \eqref{eq:long-orbit} and \eqref{eq:t-ineq} we get $\mu p_n = \sum_0^n\mu t_i > \sum_0^n\|x_{i+1}-x_i\| > \sigma$, that is, $p_n > \sigma/\mu>1-\varepsilon$ from \eqref{eq:mu-def}.

Since $p_{i+1}-p_i=t_i< \varepsilon$, there is $m\in \{0,\ldots,n\}$ such that
$$
   1 - \varepsilon < p_m < 1.
$$
Using \eqref{eq:t-ineq} and the triangle inequality we can estimate $\|x_{m+1}\|\le \sum_{i=0}^m\|x_{i+1}-x_i\| <\mu p_m <\sigma$.
That is,
\begin{equation}\label{eq:x:m+1}
    x_{m+1}\in \sigma B_X.
\end{equation}
On the other hand, \eqref{eq:y:smaller} for $i=m$ gives
\begin{eqnarray*}
  | g(x_{m+1})-\ol y|_k &\le& |g(x_{m+1})-p_m\ol y|_k+|p_m\ol y-\ol y|_k\\
  &<&  \varepsilon p_m + (1-p_m)|\ol y|_k\\
  &<& \varepsilon + \varepsilon |\ol y|_k.
\end{eqnarray*}
Recalling the choice of $\varepsilon$, we get
$$
   | g(x_{m+1})-\ol y|_k  < \varepsilon_1.
$$
This and \eqref{eq:x:m+1} mean that there exists $x\in \sigma B_X$ such that $| g(x)-\ol y|_k  <  \varepsilon_1$.

Since $\varepsilon_1>0$ and $k\ge k_0$ were arbitrary, we have that
$$
\forall \varepsilon_1 >0,   \forall k\ge k_0,\  \exists x\in \sigma B_X:\
   | g(x)-\ol y|_k  <  \varepsilon_1.
$$

Having in mind that $|\cdot|_k$ seminorms are graded, the latter implies that
$$
\forall \varepsilon_1 >0,   \forall k\ge 0,\ \exists x\in \sigma B_X:\
   | g(x)-\ol y|_k  <  \varepsilon_1.
$$

Lemma~\ref{lem1} completes the proof.
\qed

\begin{thm}\label{thm:kornerstone}
  Let $(X,\|\cdot\|_n)$ and $(Y,\|\cdot\|_n)$ be Fr\'echet spaces and let
  $$
    f:X\to Y
  $$
  be a continuous function.

  Let  $f$ be directionally differentiable on the open set $V\subset X$.

  Assume that for some $\mathbf{s}\in \mathbb{R}_+^\infty$ and some non-empty set $C\subset Y$ we have
 \begin{equation}\label{eq:main:th2:k2}
    f'(x,\Pi_{\mathbf{s}}(X))\supset  C,\quad \forall x\in V .
   \end{equation}
  Then for  any $x$ such that $x+\Pi_{\mathbf{s}}(X)\subset V$, it holds that
  \begin{equation}\label{eq:main:th2:k3}
    \cl{f(x+\Pi_{\mathbf{s}}(X))}\supset f(x) + C.
  \end{equation}
\end{thm}

\textbf{Proof.} 
 Fix $x_0\in V$ such that $x_0+  \Pi_{\s}(X){\subset }V$.

  Obviously, it is enough to prove that
  \begin{equation}\label{eq:korn:1}
    f(x_0)+\ol y\in\cl{f(x_0+\Pi_{\s}(X))}
  \end{equation}
  for each fixed $\ol y \in C$.

  We will apply Proposition~\ref{prop:4} to the Banach space $(X_{\s},\|\cdot\|_{\s})$, the Fr\'echet space $(Y,\|\cdot\|_n)$ and the function $g:X_\s\to  Y$ defined by
  $$
    g(x) := f(x_0 + x)-f(x_0), \quad\forall x\in X_\s
  $$
  while $\ol y$ will play the same role.

  To this end we will check the other assumptions of Proposition~\ref{prop:4}. In our case $\sigma$ will be equal to  $1$ and the set $U$ will be $U:=(V-x_0)\cap X_{\s}$. It is clear that $U$ is open in $X_{\s}$ (since $\|\cdot\|_{\s}$-topology is stronger than $\rho_{(X,\|\cdot\|)}$-topology). Since $x_0+  \Pi_{\s}(X)\subset V$, we also have $U\supset  \Pi_{\s}(X) = B_{X_\s}$.

  Since $B_{X_{\s}}=\Pi_{\s}(X)$, we have by definition
  $$
    g'(x,B_{X_{\s}}) = f'(x_0+x,\Pi_{\s}(X)),\quad \forall x \in X_\s.
  $$

  If $x\in U$ then  $x_0+x\in V$ and we have by \eqref{eq:main:th2:k2} that $f'(x_0+x,\Pi_{\s}(X)) \supset C$. Since $\ol y\in C$ we have, therefore,
  $$
   \ol y\in g'(x, B_{X_{\s}}),\quad\forall x\in U.
  $$
 Proposition~\ref{prop:4} implies that
  $$
  \ol y\in \cl{g( B_{X_{\s}})} = \cl{f(x_0+\Pi_\s(X))-f(x_0)},
  $$
  which is \eqref{eq:korn:1} and the proof is completed.
\qed

At this point we have all we need  in order to get the promised surjectivity result.

\begin{thm}[Local surjectivity]\label{thm:surj}
  Let  $(X,\|\cdot\|_n)$ and $(Y,\|\cdot\|_n)$ be Fr\'echet spaces and let $U\subset X$  be nonempty and open.

  Let the function $f:X\to Y$ be continuous, directionally differentiable on $U$, and  such that
  \begin{equation}
   \label{eq:kapa-inc}
   f'(x,\Pi_{\mathbf{s}}(X))\supset \Pi_{\mathbf{s}}(Y)
  \end{equation}
  for all $x\in U$, and all $\s\in \mathbb{R}_+^\infty$. Then
  \begin{equation}
   \label{eq:surj}
     f (B(x,r))\supset B^\circ(f(x),r),\quad
  \forall x\in U,\ \forall r:\ 0<r< m_U(x),
  \end{equation}
  where
  $
  m_U(x) := \mathrm{dist}(x,X\setminus U).
  $

In particular, $f$ is open at linear rate and locally surjective on $U$.
\end{thm}

\textbf{Proof.} 
If we fix $\mathbf{s}\in \mathbb{R}_+^\infty$ and set $C:= \Pi_{\mathbf{s}}(X)$, from Theorem~\ref{thm:kornerstone} we will get that for any $x$ such that $x+\Pi_{\mathbf{s}}(X)\subset U$,
  \[
  \cl f(x+\Pi_{\mathbf{s}}(X))\supset f(x)+ \Pi_{ \mathbf{s}}(Y),
  \]
  which is \eqref{eq:oct_2} and the claim follows from Proposition~\ref{pro:oct_1}.
\qed

We complete this section by rewriting the previous result in more analytical terms.

 \begin{corollary}\label{cor:principal}
 Let $(X,\|\cdot\|_n)$ be a Fr\'echet space and let $(Y,\|\cdot\|_n)$ be a non-exotic Fr\'echet space and let
  $$
    f:X\to Y
  $$
  be a continuous function. Let $U\subset X$ be open and let $f$ be directionally differentiable on $U$. Assume that there are $d\in \N\cup\{0\}$ and $c_n>0$ such that for any $x\in U$ and any $v\in Y$,
  \begin{equation}\label{eq:oh}
  \exists h\in X\ :\ f'(x,h)=v\ \mbox{and}\ \| h\|_n\le c_n\|v\|_{n+d},\ \forall n\in \N\cup\{0\}.
  \end{equation}
  Then $f$ is locally surjective on $U$.

  More precisely, if $\rho_{(Y,c_n\|\cdot\|_{n+d})}$ is the metric on $Y$ given by \eqref{eq:rho} for the equivalent system of seminorms $(c_n\|\cdot\|_{n+d})_{n\ge0}$, then  \eqref{eq:surj} holds.
  \end{corollary}

\textbf{Proof.} 
  In $Y$ define the seminorms
  $$
    |||\cdot|||_n:=c_n\|\cdot\|_{n+d},\quad n\ge 0.
  $$

  Since $Y$ is non-exotic, the metric $\rho_{(Y,\|\cdot\|_n)}$ is equivalent to  the metric $\rho_{(Y,|||\cdot\|||_n)}$.

 For $(X,\|\cdot\|_n)$ and $(Y,|||\cdot|||_n)$ condition (\ref{eq:oh}) can be rewritten as
  \[
 f'(x,\Pi_\s(X,\|\cdot\|_n))\supset \Pi_{\s}(Y,|||\cdot|||_n),\quad \forall x\in U,\   \forall \s \in \R^\infty _+.
  \]

  Indeed, if $v\in \Pi_{\s}(Y)$ then $|||v|||_n\le s_n$.

  From (\ref{eq:oh}) it follows that there is $h\in X$ such that $f'(x,h)=v$ and $\| h\| _n\le c_n\|v\|_{n+d}$. But  $c_n\|v\|_{n+d} =|||v|||_n$, hence $\| h\| _n\le   |||v|||_n\le s_n$.

  The conditions of Theorem \ref{thm:surj} hold and this completes the proof.
\qed

  Let us note that in the above corollary the assumption on $Y$ to be non-exotic is linked only with the loss of derivatives, i.e. in the case $d=0$ this assumption is redundant.

\section{Proof of the main result}\label{sec:main}

Here we prove Theorem~\ref{thm:main-main}.  We will derive it from the following partial case.

\begin{proposition}
\label{pro:main-partial}
 Let $(X,\|\cdot\|_n)$ and $(Y,\|\cdot\|_n)$ be Fr\'echet spaces and let
  $
    f:X\to Y
  $
  be a continuous function such that $f(0)=0$. Let  $U\ni 0$ be a nonempty  open subset of $X$ such that $f$ is injective and  directionally differentiable on $U$; and such that for any $x\in U$ and any $v\in Y$,
  \begin{equation}\label{eq:oh-prim}
  \exists h\in X\ :\ f'(x,h)=v\ \mbox{and}\ \| h\|_n\le \|v\|_{n},\ \forall n\ge 0.
  \end{equation}
  Then there exists an open $V\ni 0$ such that $f$ is invertible on $V$ and
  \begin{equation}
     \label{eq:invert-prim}
     \|f^{-1}(v)-f^{-1}(u)\|_n\le \|v-u\|_{n},\quad \forall u,v\in f(V),\ \forall n\ge 0.
  \end{equation}
   \end{proposition}

\textbf{Proof.} 
      Obviously \eqref{eq:oh-prim} implies that
$$
        f'(x,\Pi_\s(X))\supset \Pi_{\s}(Y),\quad \forall x\in U, \ \forall  \s \in \R^\infty _+.
   $$

Then we know from Theorem~\ref{thm:surj} that $f$ on $U$  is an open mapping (that is, $f(U_1)$ is open for each open $U_1\subset U$) and it is locally surjective on $U$. Hence, there exists a neighbourhood $V_1\subset U$ of $0$ such that $f$ is surjective and, therefore, invertible, on $V_1$. The set $f(V_1)$ is open and $f^{-1}$ is continuous on $f(V_1)$. Let $W$ be open convex set such that $f(0)\in W$ and $W\subset f(V_1)$. Set $V:=f^{-1}(W)$. Obviously $V$ is open, $V\subset U$ and $f(V) \equiv W$.

Applying Theorem~\ref{thm:kornerstone} with $C=\Pi_{\s}(Y)$ we get that
\begin{equation}
    \label{eq:prop-star}
       \cl{f(x+\Pi_{\mathbf{s}}(X))}\supset f(x) + \Pi_{\s}(Y),\quad \forall x:\ x+\Pi_{\mathbf{s}}(X) \subset V.
   \end{equation}

 It is clear that if $x+\Pi_\s(X)\subset V$, then $f(x+\Pi_\s(X))\subset f(V)$.

 We will prove the following

 \bigskip

\noindent
 \textsc{Claim.} If $\s' \in  \R^\infty _+$ is such that $x+\Pi_{\s'}(X)\subset V$ and $\cl f(x+\Pi_{\s'}(X))\subset f(V)$, then
   \[
   \cl f(x+\Pi_{\s'}(X))= f(x+\Pi_{\s'}(X)).
   \]

  Take $y\in \cl f(x+\Pi_{\s'}(X))$ and let $y=\lim y_n$ with $y_n=f(x_n)$, $x_n\in x+ \Pi_{\s'}(X)$.

  By the presumed inclusion, there exists $z\in V$ such that $y=f(z)$.

  Since $f^{-1}$ is continuous on $f(V)$, and in particular at $y$,
  \[
  x_n=f^{-1}(y_n)\to f^{-1}(y)=z.
  \]

  Since $x_n\in x+ \Pi_{\s'}(X)$, and the latter set is closed, we get that $z\in  x+ \Pi_{\s'}(X)$. That is, $y\in f( x+ \Pi_{\s'}(X))$ and the claim follows.

  \bigskip

Now let $ x+ \Pi_{\s}(X)\subset V$. Because $f$ is continuous at $x$, for any $\eps >0$ there exists $\gamma =\gamma(x)>0$ such that $B(x,\gamma)\subset V$ and
\begin{equation}\label{eq:inc}
f\left( B(x,\gamma)\right) \subset B(f(x),\eps).
\end{equation}
 We claim that there exists $\delta=\delta(x,\s,\eps)>0$ such that for $t\in (0,\delta)$,
 \[
 x+t\Pi_\s(X)\subset B(x,\gamma).
 \]

 Indeed,
$\rho_{(X,\|\cdot\|_n)}( x,x+ t\Pi_{\s}(X))=\rho_{(X,\|\cdot\|_n)}(0,t\Pi_\s(X))$ by the shift invariance of the metric $\rho_{(X,\|\cdot\|_n)}$ and $\displaystyle \limsup _{t\downarrow 0} \rho_{(X,\|\cdot\|_n)}(0, t\Pi_\s(X))=0$ as it is shown in \cite[Lemma 3.3]{NME-JOTA}.

So, we fix $\eps >0$ such that $B(f(x),\eps)\subset f(V)$. Then by \eqref{eq:inc} we get
\[
\cl f(x+t\Pi_\s(X))\subset \cl f(B(x,\gamma))\subset B(f(x),\eps)\subset f(V),
\]
for all $t\in (0,\delta)$.

Having in mind that $t\Pi_\s(X)=\Pi_{t\s}(X)$ for all $t\ge 0$ and all $\s \in \R^\infty _+$, we apply the claim for $x$ and $\s'=t\s$ to obtain that
\[
\cl f(x+t\Pi_\s(X))=f(x+t\Pi_\s(X)), \quad \forall t\in (0,\delta).
\]

 From \eqref{eq:prop-star} we get
      $$
         f(x+t\Pi_{\mathbf{s}}(X))\supset f(x) + t \Pi_{\mathbf{s}}(Y), \quad \forall t\in (0,\delta).
$$

From surjectivity of $f$, for any $y\in f(V)$ there exists $x\in V$ such that $f(x)=y$. We have actually proved that for any $\s \in  \R^\infty _+$ there exists $\delta=\delta(x(y),\s,\eps)>0$ such that
 $$
         f(x+t\Pi_{\mathbf{s}}(X))\supset y + t \Pi_{\mathbf{s}}(Y), \quad \forall t\in (0,\delta),
$$
 or, equivalently,
\begin{equation}\label{eq:star-inverse}
        f^{-1}(y+ t \Pi_{\mathbf{s}}(Y)) \subset f^{-1}(y) + t \Pi_{\mathbf{s}}(X), \quad \forall t\in (0,\delta).
      \end{equation}
We will now derive from the latter the Lipschitz-like inequalities of the type we want, but apparently weaker.

Take any  $y\in f(V)$ and any $z\in Y$, $z\neq 0$. Let $\s:= (\|z\|_n)_{n=0}^\infty$. We invoke \eqref{eq:star-inverse} with this $\s$. For  $t\in (0,\delta)$,  immediately $tz\in t\Pi_\s(Y)$ and, therefore,
\[
f^{-1}(y+ t z) - f^{-1}(y) \subset t \Pi_{\mathbf{s}}(X).
\]
Hence,
      \begin{equation}
        \label{eq:star-inequality}
        \|f^{-1}(y+ tz) -  f^{-1}(y)\|_n \le  t \|z\|_n, \quad \forall n\ge 0, \quad \forall t\in (0,\delta).
      \end{equation}

     Now, fix $n\ge 0$ and $u,v\in f(V)$, $u\neq v$. Set  $z:=v-u$. Define for $\lambda\in [0,1]$ the function
     $$
         g(\lambda) := \|f^{-1}((1-\lambda)u+\lambda v) - f^{-1}(u)\|_n.
     $$
     For each fixed $\lambda \in (0,1)$ the estimate \eqref{eq:star-inequality} applied at $y:= (1-\lambda)u+\lambda v\in f(V)$ (since $f(V)$ is a  convex set) gives for all $t\in (0,\delta)$  that
\begin{eqnarray*}
g(\lambda+t)- g(t) &=& \|f^{-1}(y+tz)-f^{-1}(u)\|_n-\|f^{-1}(y)-f^{-1}(u)\|_n\\
  &\le &\|f^{-1}(y+tz)-f^{-1}(y)\|_n \\
  &\le &t\|z\|_n.
\end{eqnarray*}
 So,
     $$
        g^+(\lambda)\le \|z\|_n.
     $$
     Lemma~\ref{lem:dini-mvt} ensures that $g(1)\le \|z\|_n$ and we are done.
\qed

 We   now have all we need to proof our  main result.

\textbf{Proof of Theorem~\ref{thm:main-main}. }
 Fix $x\in U$. By considering instead of $f$ the function $f(\cdot-x)-f(x)$ we may assume that $x=0$ and $f(0)=0$, as this will simplify the considerations.

 Since $(Y,\|\cdot\|_n)$ is non-exotic, we may consider $Y$ with the seminorms
 $$
  |||\cdot|||_n:=c_n\|\cdot\|_{n+d},\quad n\ge 0,
 $$
 as we already did in the proof of Corollary~\ref{cor:principal}.

 Obviously \eqref{oh} is equivalent to \eqref{eq:oh-prim} for the seminorms $ |||\cdot|||_n$ and, therefore, Proposition~\ref{pro:main-partial} gives \eqref{eq:invert-prim} for $ |||\cdot|||_n$ which in its turn is equivalent to \eqref{eq:invert}.
\qed

\section{Local injectivity}\label{sec:injectivity}

Here we give sufficient conditions for injectivity in the spirit of \cite{hamilton}. If we were to assume tameness of the space, we could derive these from seemingly more natural conditions, like in \cite{hamilton}, but we would rather leave this to the readers who are specifically interested in the context set by \cite{hamilton}.

We would be more interested in general sufficient conditions for injectivity in therms of directional derivatives in a Banach space, as a first step towards the general Fr\'echet case; but this might be a topic of further research.

 \begin{thm}[Local injectivity]\label{injectivity}
 Let $(X,\|\cdot\|_n)$ and $Y(\|\cdot\|_n)$ be  Fr\'echet spaces,  $f:U\subset X\to Y$ be a continuous function such that $f\in C^1(U)$, and $U$ be an open and convex set.

If there exist some constants $c_n>0$ such that
\begin{equation}\label{eqq1}
\| h\| _n\le c_n\| f'(x)(h)\| _{n},\quad \forall n\ge 0, \forall x\in U,\ \forall h\in X,
\end{equation}
and there exist some constants $c_n'>0$ such that
\begin{equation}\label{eqq2}
\begin{array}{l}
\| f'(x)(h)-f'(z)(h)\| _n\le c_n'[ \| x-z\| _{r}\|h\|_n+\| x-z\| _{n}\|h\|_r],\\[10pt]
\forall n\ge 0,\ \forall x,z\in U,\ \forall h\in X, \mbox{ and some } r \ge 0,
\end{array}
\end{equation}
then $f$ is locally injective on $U$.
\end{thm}

\textbf{Proof. }
Fix $\ol x\in U$. Let $U'$ be a neighbourhood of $\ol x$ contained in $U$.

For $x_0,x_1\in U'$, from Taylor's theorem with integral reminder,
\[
f(x_1)=f(x_0)+\int_0^1 f'(x_0+t(x_1-x_0))(x_1-x_0)\,dt.
\]

Adding and subtracting  $f'(x_0)(x_1-x_0)$ we get
\[
f(x_1)=f(x_0)+f'(x_0)(x_1-x_0)+\int_0^1 [f'(x_0+t(x_1-x_0))-f'(x_0)](x_1-x_0)\,dt.
\]

Let $x_0,x_1\in U'$. Then,
\[
f'(x_0)(x_0-x_1)=f(x_0)-f(x_1)+\int_0^1 [f'(x_0+t(x_1-x_0))-f'(x_0)](x_1-x_0)\,dt,
\]
and for any $n\ge 0$,
\[
\|f'(x_0)(x_1{-}x_0)\|_n{=}\left \| f(x_0){-}f(x_1){+}\int_0^1 [f'(x_0{+}t(x_1{-}x_0)){-}f'(x_0)](x_1{-}x_0)\,dt\right \|_n .
\]
Hence, using \eqref{eqq1},
\begin{equation}\label{eqq3}
\begin{array}{rcl}
\displaystyle \frac{1}{c_n}\| x_0-x_1\| _n&\le &\|f'(x_0)(x_1-x_0)\|_n\\
 &\le &  \| f(x_0)-f(x_1)\|_n{+}\\
 &  & \displaystyle \int
_0^1 \|[f'(x_0+t(x_1-x_0))-f'(x_0)](x_1-x_0)\|_n \,dt \\
  &\le &\displaystyle \| f(x_0)-f(x_1)\|_n +\\
  & &\displaystyle \int_0^1 2c_n' t \| x_1-x_0  \|_r\| x_1-x_0\|_n \,dt\\
 &\le & \| f(x_0)-f(x_1)\|_n +c_n'  \| x_1-x_0  \|_r\| x_1-x_0\|_n.\\
\end{array}
\end{equation}

Let $U'':=U'\cap \{ x:\| x-\ol x\| _r< \delta \}$, where $\displaystyle \delta <\frac{1}{2c_rc_r'}$. For $x_0,x_1\in U''$ from \eqref{eqq3} taken for $n=r$  we get
 \[
 \begin{array}{rcl}
\displaystyle  \frac{1}{c_r}\| x_0-x_1\| _r&\le &\| f(x_0)-f(x_1)\|_r+c_r'\| x_1-x_0\|_r^2\\[10pt]
  &\le & \| f(x_0)-f(x_1)\|_r+c_r'[\| x_1-\ol x\|_r+\|x_0-\ol x\|_r]\| x_1-x_0\|_r\\[10pt]
  &\le &  \| f(x_0)-f(x_1)\|_r +2c_r'\delta  \| x_1-x_0\|_r,\\
  \end{array}
  \]
 or
 \[
 \left(\frac{1}{c_r}-2\delta c_r'\right)\| x_0-x_1\| _r\le \| f(x_0)-f(x_1)\|_r.
 \]
 By the choice of $\delta$, for some positive $c>0$,
  \begin{equation}\label{eqq4}
 \| x_0-x_1\| _r\le c\| f(x_0)-f(x_1)\|_r.
 \end{equation}

Finally, using the estimate \eqref{eqq4} in \eqref{eqq3} we obtain that
\[
\| x_0-x_1\| _n \le c_n\| f(x_0)-f(x_1)\|_n+c_nc_n' c\| f(x_0)-f(x_1)\|_r\| x_1-x_0\|_n.
\]
The latter yields that if $f(x_0)=f(x_1)$ then $x_1=x_0$, hence $f$ is injective on~$U''$.
\qed

\section{Implicit Map Theorem}\label{sec:ift}

As well known, Inverse and Implicit Function type theorems are closely related. Here we follow the method of \cite{ngth} to derive from Corollary~\ref{cor:principal} an Implicit Function Theorem or,  more precisely, \textit{Implicit Map Theorem}, because we do not get an unique solution in general. Of course, our result is slightly different than that of \cite{ngth}, because the Inverse Function Theorem we use as a base is different.

\begin{thm}\label{thm:ift}
    Let $(X,\|\cdot\|_n)$ and $(Y,\|\cdot\|_n)$ be Fr\'echet spaces such that $Y$ is non-exotic, and let $P$ be a topological space. Consider a function
    $$
        f: X\times P\to Y,\mbox{ such that }f(\bar x, \bar p) = 0,
    $$
    where $(\bar x, \bar p)\in X\times P$ is fixed.

    Assume that there is a neighourhood $U$ of $\bar x$ in $X$ such that for any fixed  parameter $p\in P$ the function $f(\cdot, p): X\to Y$ is continuous and directionally differentiable on $U$. Denote the directional derivative of this function at $x\in U$ in the direction $h\in X$ by $f'((x,p),h)$.

    Assume also that:

    $(\i)$ $f$ is continuous at $(\bar x, \bar p)$, and

    $(\i\i)$ there are $c_n>0$ and $d\ge 0$ such that for all $x\in U$, $p\in P$ and $v\in Y$, there is $h\in X$ such that
    \begin{equation}
     \label{eq:double-triangle}
       f'((x,p),h) =v\mbox{ and }\|h\|_n\le c_n\|v\|_{n+d},\ \forall n\ge 0.
    \end{equation}
    Then there are a neighbourhood $U'$ of $\bar x$ in $X$ and a neighbourhood $O$ of $p$ in $P$ such that
    \begin{equation}
     \label{eq:ift-conclusion}
     d_{\rho_{(X,\|\cdot\|_n)}}(x,S(p)) \le \rho_{(Y,c_n\|\cdot\|_{n+d})}(0,f(x,p)),\quad \forall (x,p)\in U'\times O,
    \end{equation}
    where
    $$
        S(p) := \{ x\in X:\ f(x,p) = 0\},
    $$
    and $\rho_{(Y,c_n\|\cdot\|_{n+d})}$ is the metric on $Y$ given by \eqref{eq:rho} for the equivalent system of seminorms $(c_n\|\cdot\|_{n+d})_{n\ge0}$.

    In particular, $S(p)$ is not empty for $p\in O$.
\end{thm}

\textbf{Proof. }
   Let $\varepsilon>0$ be such that
   $$
        B(\bar x,2\varepsilon) \subset U.
   $$
   Since $f$ is equal to zero and continuous at $(\bar x, \bar p)$, there are $\delta\in(0,\varepsilon)$ and an open neighbourhood $O$ of $\bar p$ such that
   $$
        f(B(\bar x,\delta), O) \subset B^\circ (0,\varepsilon).
   $$
   Set
   $$
        U' := B^\circ (\bar x,\delta).
   $$
   Fix arbitrary $x\in U'$, and $p\in O$. Since $\rho_{(Y,c_n\|\cdot\|_{n+d})}(0,f(x,p)) <  \varepsilon$, there is  $r>0$ such that
   \begin{equation}\label{eq:ift}
        \rho_{(Y,c_n\|\cdot\|_{n+d})}(0,f(x,p)) < r < \varepsilon.
   \end{equation}
  Since $B(x,r)\subset B(x,\eps)\subset B(\ol x, 2\eps)\subset U$, Corollary~\ref{cor:principal} and \eqref{eq:ift} give that \eqref{eq:surj} is satisfied for $f(\cdot,p)$, that is,
  $$
        f(B(x,r),p)\supset B^\circ (f(x,p),r).
  $$
  From \eqref{eq:ift} it follows that $0\in B^\circ(f(x,p),r)$, so $0\in f(B(x,r),p)$, that is,
  $$
        d_{\rho_{(X,\|\cdot\|_{n})}}(x,S(p)) \le r
  $$
  and \eqref{eq:ift-conclusion} follows.
\qed

\section{Application to differential equations}\label{sec:ode}

Here we generalise Theorem~8 from \cite{ngth} by dropping the assumption that the spaces should be \textit{standard}. In essence we just put Theorem~\ref{thm:ift} in the place of the Implicit Function Theorem used in \cite{ngth}, so we do not need to write the proof in much detail.

If $(X,\|\cdot\|_n)$ is a Fr\'echet space then the space of continuous functions from $[-1,1]$ into $X$, denoted by $C([-1,1],X)$, is also Fr\'echet with the seminorms
$$
    \|u(\cdot)\| _n = \sup_{t\in[-1,1]} \|u(t)\| _n.
$$
$C^1$ is defined in a similar fashion.
Obviously, if $X$ is graded, so are  $C([-1,1],X)$ and $C^1([-1,1],X)$.

Where needed for clarity, we will indicate the space to which certain seminorm belongs by superscript, for example
$$
   \|u(\cdot)\| _n^{C^1([-1,1],X)} = \sup_{t\in[-1,1]}\max\{\|u(t)\| _n, \|u'(t)\| _n\}.
$$

\begin{thm}\label{thm:ode}
   Let $U$ be an open subset of the graded Fr\'echet space $(X,\|\cdot\|_n)$ and let $f:\mathbb{R}\times U\to X$ be continuous on $[-r_0,r_0]\times U$.

   Assume that for all $t\in[-r_0,r_0]$ the function $f(t,\cdot)$ is G\^ateaux differentiable on $U$ with derivative $D_xf(t,\cdot)$ and for some $c_n\ge 0$
   \begin{equation}
    \label{eq:ode-condition}
     \| D_x f(t,x)(h)\|_n \le c_n \|h\|_n,\quad\forall |t|\le r_0,\ x\in U,\ h\in X,\ n\ge 0.
   \end{equation}
   Then there is $r\in (0,r_0]$ such that the Cauchy problem
   \begin{equation}\label{eq:Cauchy}
        \begin{cases}
           x'(t) = f(t,x(t)),\quad |t|\le r,\\
           x(0) = x_0.
        \end{cases}
   \end{equation}
  has a solution $x(\cdot)\in C^1([-r,r],X)$.
\end{thm}

\textbf{Proof.} 
 By changing the variable $x(t)\to (x(t)-x_0)$ we can assume without loss of generality that $x_0=0$ and we will do so for the sake of simplicity. Further, using $t=rs$ and $z(s)=x(rs)$, $|s|\le 1$ we rewrite \eqref{eq:Cauchy} as
 \begin{equation}\label{eq:Cauchy_2}
        \begin{cases}
           z'(s) = rf(rs,z(s)),\quad |s|\le 1,\\
           z(0) = 0.
        \end{cases}
   \end{equation}

 Let $W$ be the open subset of the graded Fr\'echet space $C^1([-1,1],X)\times\mathbb{R}$ defined by
 $$
    W := \{ (z,r):\ |r|< |r_0|,\ z(s)\in U,\ \forall s\in [-1,1]\}.
 $$
 Let $F:W\to C([-1,1],X)\times X$ be
 $$
    F(z,r) := (z'(s) - rf(rs,z(s)),z(0)).
 $$
 Obviously, $F(0,0)=(0,0)$. Moreover, $z\in C^1([-1,1],X)$ is a solution of \eqref{eq:Cauchy_2} for some  $r\in (0,r_0)$ exactly when  $F(z,r)=(0,0)$.

 We will apply Theorem~\ref{thm:ift}  to show that $F(z,r)=0$ has solutions for small enough $|r|$.

 For each fixed $(z,r)\in W$ the function $F(\cdot, r)$ is G\^ateaux differentiable at $z$ with derivative
 $$
    D_z F(z,r)(u) = (u'(s) - rD_x f(rs,z(s))u(s), u(0)),\quad |s|\le 1,
 $$
 where $u\in  C^1([-1,1],X)$.

 So, for any fixed right hand side $(v,v_0)\in C([-1,1],X)\times X$, the equation
 \begin{equation}
  \label{eq:ode-inverse}
 D_z F(z,r)(u) =  (v,v_0)
 \end{equation}
 is linear ODE with respect to  $u$ with continuous linear operator thanks to \eqref{eq:ode-condition}, so from  \cite[Proposition 3.4]{popp} it follows that it has unique solution $u\in C^1([-1,1],X)$.

 Then Gronwall Lemma together with  \eqref{eq:ode-condition} is used as in the proof of Theorem~8 of \cite{ngth} to show that
 $$
    \|u\|_n^{C^1([-1,1],X)} \le (1+(1+r_0c_n)e^{r_0c_n})\|(v,v_0)\|_n^{C([-1,1],X)\times X},\quad\forall n\ge 0.
 $$

 This means that \eqref{eq:double-triangle} is fulfilled with $d=0$ and we can apply Theorem~\ref{thm:ift} to complete the proof.
\qed

\end{document}